\documentclass{amsart}

\usepackage{amsfonts,amssymb,verbatim,amsmath,amsthm,latexsym,textcomp,amscd}
\usepackage{latexsym,amsfonts,amssymb,epsfig,verbatim}
\usepackage{amsmath,amsthm,amssymb,latexsym,graphics,textcomp}
\usepackage{graphicx}
\usepackage{color}
\usepackage{url}

\begin{document}

\newtheorem{theorem}{Theorem}[section]
\newtheorem{prop}[theorem]{Proposition}
\newtheorem{lemma}[theorem]{Lemma}
\newtheorem{cor}[theorem]{Corollary}
\newtheorem{defn}[theorem]{Definition}
\newtheorem{conj}[theorem]{Conjecture}
\newtheorem{claim}[theorem]{Claim}
\newtheorem{example}[theorem]{Example}
\newtheorem{rem}[theorem]{Remark}
\newtheorem{rmk}[theorem]{Remark}
\newtheorem{obs}[theorem]{Observation}
\newtheorem{qn}[theorem]{Question}

\newcommand{\map}{\rightarrow}
\newcommand{\C}{\mathcal C}
\newcommand\AAA{{\mathcal A}}
\newcommand\BB{{\mathcal B}}
\newcommand\DD{{\mathcal D}}
\newcommand\EE{{\mathcal E}}
\newcommand\FF{{\mathcal F}}
\newcommand\GG{{\mathcal G}}
\newcommand\HH{{\mathcal H}}
\newcommand\I{{\stackrel{\rightarrow}{i}}}
\newcommand\J{{\stackrel{\rightarrow}{j}}}
\newcommand\K{{\stackrel{\rightarrow}{k}}}
\newcommand\LL{{\mathcal L}}
\newcommand\MM{{\mathcal M}}
\newcommand\NN{{\mathbb N}}
\newcommand\OO{{\mathcal O}}
\newcommand\PP{{\mathcal P}}
\newcommand\QQ{{\mathcal Q}}
\newcommand\RR{{\mathcal R}}
\newcommand\SSS{{\mathcal S}}
\newcommand\TT{{\mathcal T}}
\newcommand\UU{{\mathcal U}}
\newcommand\VV{{\mathcal V}}
\newcommand\WW{{\mathcal W}}
\newcommand\XX{{\mathcal X}}
\newcommand\YY{{\mathcal Y}}
\newcommand\ZZ{{\mathbb Z}}
\newcommand\hhat{\widehat}
\newcommand\vfn{\stackrel{\A}{r}(t)}
\newcommand\dervf{\frac{d\stackrel{\A}{r}}{dt}}
\newcommand\der{\frac{d}{dt}}
\newcommand\vfncomp{f(t)\I+g(t)\J+h(t)\K}
\newcommand\ds{\sqrt{f^{'}(t)^2+g^{'}(t)^2+h^{'}(t)^2}dt}
\newcommand\rvec{\stackrel{\A}{r}}
\newcommand\velo{\frac{d\stackrel{\A}{r}}{dt}}
\newcommand\speed{|\velo|}
\newcommand\velpri{\rvec \,^{'}}

\title{ Graphs of hyperbolic groups and a limit set intersection theorem}

\author{Pranab Sardar}
\address{Indian Institute of Science Education and Research Mohali}
\thanks{2010 {\em Mathematics Subject Classification.} Primary 20F67}

\thanks{{\em Key words and phrases.} Hyperbolic groups, limit sets, Bass-Serre theory}
\date{\today}

\maketitle

\begin{abstract}
We define the notion of {\em limit set intersection property} for a collection of subgroups of a hyperbolic group; namely,
for a hyperbolic group $G$ and a collection of subgroups $\mathcal S$ we say that $\mathcal S$
satisfies the {\em limit set intersection property} if for all $H,K \in \mathcal S$ we have
$\Lambda(H)\cap \Lambda(K)=\Lambda(H\cap K)$. Given a hyperbolic group admitting a decomposition into a finite graph of
hyperbolic groups structure with QI embedded condition, we show that the set
of conjugates of all the vertex and edge groups satisfy the limit set intersection property.
\end{abstract}

\section{Introduction}
Limit set intersection theorems first appear in the work of Susskind and Swarup (\cite{suss-swarup}) in the context of
geometrically finite Kleinian groups. Later on Anderson (\cite{anderson-limset1}, \cite{anderson-limset2}) undertook
a detailed study of this for general Kleinian groups. In the context of (Gromov) hyperbolic groups this is true
for quasiconvex subgroups(see \cite{GMRS}, Lemma 2.6). (Recently W. Yang has looked at the case of relatively quasiconvex
subgroups of relatively hyperbolic groups. See \cite{yang-limset}.) However, this theorem is false for general subgroups of
hyperbolic groups and no characterizations other than quasi-convexity are known
for a pair of subgroups $H,K$ of a hyperbolic group $G$ which guarantee that $\Lambda(H)\cap \Lambda(K)=\Lambda(H\cap K)$.
This motivates us to look for subgroups other than quasiconvex subgroups which satisfy the limit set intersection property.
Our starting point is the following celebrated theorem of Bestvina and Feighn.

\begin{theorem}(\cite{BF}) \label{bf-thm}
Suppose $(\mathcal G, Y)$ is a graph of hyperbolic groups with QI embedded condition and the
hallways flare condition. Then the fundamental group, say $G$, of this graph of groups is hyperbolic.
\end{theorem}

Graphs of groups are briefly recalled in section 3. There are many examples of hyperbolic groups admitting such a decomposition
into graphs of groups where the vertex or edge groups are not quasiconvex. Nevertheless, with the terminologies of the above
theorem, we have:

{\bf Theorem.} {\em The set of all conjugates of the vertex and edge groups of $G$ satisfy the limit set intersection property.}

We note that a special case of our theorem was already known by results of Ilya Kapovich (\cite{ilya-kapovich}).
There the author showed that given a $k$-acylindrical graph of hyperbolic groups $(\mathcal G, Y)$ with quasi-isometrically
embedded condition and with fundamental group $G$ (which turnns out to be hyperbolic by Theorem \ref{bf-thm},)
the vertex groups are quasiconvex subgroups of $G$. Hence, the conjugates of all the vertex groups satisfy the limit set
intersection property in this case.

\smallskip

\noindent {\bf Acknowledgments:} This work started during the author's post doctoral term in
the University of California, Davis. The author would like to thank Michael Kapovich for many
helpful discussions in this regard. The author would also like to thank Mahan Mj for helpful discussions
and Ilya Kapovich for useful email correspondence. This research was supported by the University of California, Davis.
It is partially supported by the DST INSPIRE faculty award {\small $DST/INSPIRE/04/2014/002236$}.

\section{Boundary of Gromov hyperbolic spaces and Limit sets of subspaces}

We assume that the reader is familiar with the basics of (Gromov) hyperbolic metric spaces and the coarse language.
We shall however recall some basic definitions and results that will be explicitly used in the sections to follow. 
For details one is referred to \cite{gromov-hypgps} or \cite{bridson-haefliger}.

{\bf Notation and convention.} 
In this section we shall assume that all the hyperbolic metric spaces are proper
geodesic metric spaces. 
We use QI to mean both {\em quasi-isometry} and {\em quasi-isometric} depending on
the context. Hausdorff distance of two subsets $A,B$ of a metric space $Z$ is denoted by $Hd(A,B)$.
For any subset $A$ of a metric space $Y$ and any $D\geq 0$, $N_D(A)$ will denote the $D$-neighborhood of $A$
in $Y$.
{\em We assume that all our groups are finitely generated.} For any graph $Y$ we shall denote by $V(Y)$ and $E(Y)$
the vertex and the edge sets of $Y$ respectively.

\begin{defn}
$1.$ Suppose $G$ is a group generated by a finite set $S\subset G$ and let $\gamma\subset \Gamma(G,S)$ be a path
joining two vertices $u,v\in \Gamma(G,S)$. Let $u_0=u,u_1,u_2,...,u_n=v$ be the consecutive vertices on $\gamma$.
Let $u_{i+1}=u_ix_i$, $x_i\in S\cup S^{-1}$ for $0\leq i\leq n-1$. {\em Then we shall say that the word $w=x_0x_1...x_{n-1}$
labels the path $\gamma$}. 

$2.$ Also, given $w\in \mathbb F(S)$- the free group on $S$, its image in $G$ under the natural map $\mathbb F(S)\map G$
will be called {\em the element of $G$ represented by $w$}. 
\end{defn}

\begin{defn} (See \cite{bridson-haefliger})
$1.$ Let $X$ be a hyperbolic metric space and $x\in X$ be a base point. Then the (Gromov) boundary $\partial X$
of $X$ is the equivalence classes of geodesic rays $\alpha$ such that $\alpha(0)=x$ where
two geodesic rays $\alpha, \beta$ are said to be equivalent if $Hd(\alpha,\beta)< \infty$.

The equivalence class of a geodesic ray $\alpha$ is denoted by $\alpha(\infty)$.

$2.$ If $\{x_n\}$ is an unbounded sequence of points in $X$, we say that $\{x_n\}$ converges to some
boundary point $\xi\in \partial X$ if the following holds: Let $\alpha_n$ be any geodesic joining
$x$ to $x_n$. Then any subsequence of $\{\alpha_n\}$ contains a subsequence uniformly converging on compact sets
to a geodesic ray $\alpha$ such that $\alpha(\infty)=\xi$. In this case, we say that $\xi$ is the limit of $\{x_n\}$
and write $\lim_{n\map \infty} x_n=\xi$.

$3.$ The limit set of a subset $Y$ of $X$ is the set $\{\xi\in \partial X: \exists \, \{y_n\}\subset Y \mbox{with} \, 
\lim_{n\map \infty} y_n=\xi\}$. We denote this set by $\Lambda(Y)$.
\end{defn}
The following lemma is a basic exercise in hyperbolic geometry and so we mention it without proof. It basically
uses the thin triangle property of hyperbolic metric spaces.
(See \cite{bridson-haefliger}, Chapter III.H, Exercise 3.11.)

\begin{lemma}\label{same-limit}
Suppose $\{x_n\}, \{y_n\}$ are two sequences in a hyperbolic metric space $X$ both converging to some points of 
$\partial X$. If $\{d(x_n,y_n)\}$ is bounded then $\lim_{n\map \infty}x_n=\lim_{n \map \infty} y_n$.
\end{lemma}

\begin{lemma}\label{bdry-top}
$1.$ There is a natural topology on the boundary $\partial X$ of any proper hyperbolic metric space $X$ with respect to which
$\partial X$ becomes a compact space.

$2.$ If $f:X\map Y$ is a quasi-isometric embdedding of proper hyperbolic metric spaces then $f$ induces
a topological embedding $\partial f: \partial Y\map \partial X$.

If $f$ is a quasi-isometry then $\partial f$ is a homeomorphism.
\end{lemma}

We refer the reader to Proposition $3.7$ and Theorem $3.9$ in Chapter $III.H$ of \cite{bridson-haefliger}
for a proof of Lemma \ref{bdry-top}.

\begin{defn}
$1.$ A map $f:Y\map X$ between two metric spaces is said to be a proper embedding if for all $M>0$ there
is $N>0$ such that $d_X(f(x),f(y))\leq M$ implies $d_Y(x,y)\leq N$ for all $x,y\in Y$.

A family of proper embeddings between metric spaces $f_i:X_i\map Y_i$, $i\in I$- where $I$ is an indexing set, is said to be uniformly
proper if for all $M>0$ there is an $N>0$ such that for all $i\in I$ and $x,y\in X_i$,
$d_{Y_i}(f_i(x),f_i(y))\leq N$ implies that $d_{X_i}(x,y)\leq M$.

$2.$ If $f:Y\map X$ is a proper embeddings of proper hyperbolic metric spaces then we say that {\em Cannon-Thurston({\bf CT}) map} exists
for $f$ if $f$ gives rise to a continuous map $\partial f:\partial Y\map \partial X$.
\end{defn}
This means that given a sequence of points $\{y_n\}$ in $Y$ converging to $\xi\in \partial Y$, the sequence $\{f(y_n)\}$ 
converges to a point of $\partial X$ and the resulting map is continuous. {\em Note that our terminology is slightly different from
Mitra(\cite{mitra-trees})}. The following lemma is immediate.

\begin{lemma}\label{CT-limset}
Suppose $X,Y$ are hyperbolic metric spaces and $f:Y\map X$ is a proper embedding.
If the CT map exists for $f$ then we have $\Lambda(f(Y))=\partial f(\partial Y)$.
\end{lemma}

We mention the following lemma with brief remarks about proofs since it states some standard facts from hyperbolic geometry.

\begin{lemma}\label{bdry-lemma}
Suppose $Z$ is proper $\delta$-hyperbolic metric space and $\{x_n\}$ and $\{y_n\}$ are two sequences in $Z$
such that $\lim_{n\map\infty} x_n=\lim_{n\map \infty}y_n=\xi\in\partial X$. For each $n$ let $\alpha_n, \beta_n$
be two geodesics in $X$ joining $x_1,x_n$ and $y_1,y_n$ respectively.

$(1)$ Then there is subsequences $\{n_k\}$ of natural numbers such that the sequences
of geodesics $\{\alpha_{n_k}\}$ and $\{\beta_{n_k}\}$ converge uniformly on compact sets to two geodesics $\alpha,\beta$
joining $x_1,y_1$ respectively to $\xi$. 

$(2)$ Moreover, there is a constant $D$ depending only on $\delta$ and $d(x_1,y_1)$
and  sequences of points $p_{n_k}\in \alpha_{n_k}$ and $q_{n_k}\in \beta_{n_k}$ such that $d(p_{n_k}, q_{n_k})\leq D$,
and $\lim_{k\map\infty} p_{n_k}=\lim_{k\map \infty}q_{n_k}=\xi$.

$(3)$ The conclusion $(2)$ remains valid if we replace $\alpha_n,\beta_n$ by $K$-quasi-geodesics for some $K\geq 1$; in other words
if $x_n,y_n$ are joined to $x_1,y_1$ by $K$-quasi-geodesics $\alpha_n,\beta_n$ respectively then
there is a constant $D$ depending on $\delta$, $d(x_1,y_1)$ and  $K$, a subsequence $\{n_k\}$ of natural numbers and
sequences of points $p_{n_k}\in \alpha_{n_k}$, $q_{n_k}\in \beta_{n_k}$ such that $d(p_{n_k}, q_{n_k})\leq D$,
and $\lim_{k\map\infty} p_{n_k}=\lim_{k\map \infty}q_{n_k}=\xi$.
\end{lemma}

For a proof of $(1),(2)$ see Lemma 3.3, and Lemma 3.13 and for $(3)$ Theorem 1.7(stability of quasi-geodesics)
in Chapter III.H of \cite{bridson-haefliger}. More precisely, for proving $(3)$ we may choose geodesic segments $\alpha^{'}_n$'s,
$\beta^{'}_n$'s connecting the endpoints of the quasi-geodesics $\alpha_n$'s and $\beta_n$'s respectively and then apply $(1)$
for these geodesics to extract subsequences $\{\alpha^{'}_{n_k}\}$ and $\{\beta^{'}_{n_k}\}$ of $\{\alpha^{'}_n\}$ and $\{\beta^{'}_n\}$ 
respectively both converging uniformly on compact sets.
Then we can find two sequences of points $p^{'}_{n_k}\in \alpha^{'}_{n_k}$, $q^{'}_{n_k}\in \beta^{'}_{n_k}$ satisfying $(2)$.
Finally, by stability of quasi-geodesics for all $k$ there are $p_{n_k}\in \alpha_{n_k}$, $q_{n_k}\in \beta_{n_k}$ such that
$d(p_{n_k},p^{'}_{n_k})$ and $d(q_{n_k},q^{'}_{n_k})$ are uniformly small. That will prove $(3)$.

\begin{defn}
%\begin{enumerate}
%\item
Suppose $G$ is a Gromov hyperbolic group.
Let $\mathcal S$ be any collection of subgroups of  $G$. We say that $\mathcal S$
has the {\em limit set intersection property} if for all $H, K\in \mathcal S$ we have
$\Lambda(H)\cap \Lambda(K)=\Lambda(H\cap K)$.
%\item
%\end{enumerate}
\end{defn}

We state two elementary results on limit sets for future use.

\begin{lemma}\label{int-lemma}
Suppose $G$ is a hyperbolic group and $H$ is any {\em subset} of $G$. Then for all $x\in G$ we have

$(1)$ $\Lambda(xH)=\Lambda(xHx^{-1})$.

$(2)$ $\Lambda(xH)= x\Lambda(H)$.
\end{lemma}

$Proof:$ $(1)$ follows from Lemma \ref{same-limit}. For $(2)$ one notes that $G$ acts naturally on a Cayley graph $X$ of $G$
by isometries and thus by homeomorphisms on $\partial X=\partial G$ by Lemma \ref{bdry-top}. $\Box$

\section{Graphs of groups}

We presume that the reader is familiar with the Bass-Serre theory. However, we  briefly recall all the concepts
that we shall need. For details one is referred to section 5.3 of J.P. Serre's book {\em Trees} (\cite{serre-trees}). 
Although we always work with nonoriented metric graphs like Cayley graphs, we need oriented graphs possibly
with multiple edges between adjacent vertices and loops to describe graphs of groups. Hence the following definition
is quoted from \cite{serre-trees}.

\begin{defn}
A graph $Y$ is a pair $(V,E)$ together with two maps 

$E\map V\times V$,\hspace{0.5cm} $e\mapsto (o(e),t(e))$ and

$E\map E$, \hspace{1cm}$e\mapsto \bar{e}$

such that $o(\bar{e})=t(e), t(\bar{e})=o(e)$ and $\bar{\bar{e}}=e$ for all $e\in E$.
\end{defn}

For an edge $e$ we refer to $o(e)$ as the origin and $t(e)$ as the terminus of $e$; the edge $\bar{e}$ is the same
edge $e$ with opposite orientation. We write $V(Y)$ for $V$ and $E(Y)$ for $E$. We refer to $V(Y)$ as the set of
vertices of $Y$ and $E(Y)$ as the set of edges of $Y$. We shall denote by $|e|$ the edge $e$ without any orientation. 

\begin{defn}
A {\bf graph of groups} $(\mathcal G,Y)$ consists of the following data: 

$(1)$ a (finite) graph $Y$ as defined above,

$(2)$ for all $v\in V(Y)$ (and edge $e\in E(Y)$) there is group $G_v$ (respectively $G_e$) together with
two injective homomorphisms $\phi_{e,0(e)}: G_e\map G_{o(e)}$ and $\phi_{e,t(e)}: G_e\map G_{t(e)}$ for all $e\in E(Y)$
such that the following conditions hold: 

$(i)$ $G_e=G_{\bar{e}}$,

$(ii)$ $\phi_{e,o(e)}=\phi_{\bar{e},t(\bar{e})}$ and $\phi_{e,t(e)}=\phi_{\bar{e},o(\bar{e})}$. 
\end{defn}

We shall refer to the maps $\phi_{e,v}$'s as the canonical maps of the graph of groups.
We shall refer to the groups $G_v$ and $G_e$, $v\in V(Y)$ and $e\in E(Y)$ as vertex groups and edge groups respectively.
For topological motivations of graph of groups and the following definition of the fundamental group of a graph of groups
one is referred to \cite{scott-wall} or \cite{altop-hatcher}.

\begin{defn}{\bf Fundamental group of a graph of groups}
Suppose $(\mathcal G, Y)$ is a graph of groups where $Y$ is a (finite) connected oriented graph.
Let $T\subset Y$ be a maximal tree.
Then the fundamental group $G=\pi_1(\mathcal G,Y, T)$ of $(\mathcal G,Y)$ is defined in terms of generators
and relators as follows:

The generators of $G$ are the elements of the disjoint union of the generating sets of the vertex groups $G_v$, $v\in V(Y)$
and the set $E(Y)$ of {\em oriented} edges of $Y$. 

The relators are of four of types: 
$(1)$ Those coming from the vertex groups; $(2)$ $\bar{e}=e^{-1}$ for all edge $e$, $(3)$ $e=1$ for $|e|\in T$ and $(4)$
$e\phi_{e,t(e)}(a)e^{-1}=\phi_{e,o(e)}(a)$ for all oriented edge $e$ and $a\in G_e$.

\end{defn}

{\bf Bass Serre tree of a graph of groups}

Suppose $(\mathcal G, Y)$ is a graph of groups and let $T$ be a maximal tree in $Y$ as in the above definition.
Let $G=\pi_1(\mathcal G,Y,T)$ be the fundamental group of the graph of groups. 
The Bass-Serre tree, say $\mathcal T$, is the tree with
vertex set $\bigsqcup_{v\in V(Y)} G/G_v$ and
edge set $\bigsqcup_{e\in E(Y)} G/G^e_e$ where $G^e_e=\phi_{e,t(e)}(G_e)<G_{t(e)}$. The edge relations are given by

$$t(gG^e_e)=geG_{t(e)},\,o(gG^e_e)=gG_{o(e)}....(*)$$
Note that when $|e|\in T$ then we have $e=1$ in $G$.

{\bf Tree of metric spaces from a graph of groups}
 
Given a graph of groups $(\mathcal G, Y)$ and a maximal tree $T\subset Y$ one can form in a natural way a graph, say $X$,
on which the fundamental group $G=\pi_1(\mathcal G,Y,T)$ acts by isometries properly and cocompactly and which admits a
simplicial Lipschitz $G$-equivariant map $X\map \mathcal T$. The construction of $X$ can be described as follows.

We assume that  $Y$ is a finite connected graph and all the vertex groups and the edge groups are finitely generated.
We fix a finite generating
set $S_v$ for each one of the vertex groups $G_v$; similarly for each edge groups $G_e$ we fix a finite generating set $S_e$ and
assume that $\phi_{e,t(e)}(S_e)\subset S_{t(e)}$ for all $e\in E(Y)$.
Let $S=\cup_{v\in V(Y)} S_v\cup (E(Y)\setminus E(T))$ be a generating set of $G$ where in $E(Y)\setminus E(T)$ we shall include
only nonoriented edges of $Y$ not in $T$. We define $X$ from the disjoint union of the following graphs by introducing some
extra edges as follows:

$(1)$ {\bf Vertex spaces:} For all ${\tilde v}=gG_v\in V(\mathcal T)$, where $v\in Y$ and $g\in G$ we let $X_{\tilde v}$ denote
the subgraph of $\Gamma(G,S)$ with vertex set the coset $gG_v$; two vertices $gx,gy\in X_{\tilde v}$ are connected by an
edge iff $x^{-1}y\in S_v$. We shall refer to these subspaces of $X$ as vertex spaces.

$(2)$ {\bf Edge spaces:} Similarly for any edge ${\tilde e}=gG^e_e$ of $\mathcal T$, let $X_{\tilde e}$ denote the subgraph of
$\Gamma(G,S)$ with vertex set $geG^e_e$ where where two vertices $gex,gey$ are connected by an edge iff $x^{-1}y\in \phi_{e,t(e)}(S_e)$.
We shall refer to these subspaces of $X$ as edge spaces.

$(3)$ The extra edges connect the edge spaces with the vertex spaces as follows:

For all edge ${\tilde e}=gG^e_e$ of $\mathcal T$ connecting the vertices $\tilde{u}=gG_{o(e)}$ and $\tilde{v}=geG_{t(e)}$ of
$\mathcal T$, and $x\in G^e_e$ join $gex\in X_{\tilde e}=geG^e_e$ to $gex\in X_{\tilde v}=geG_{t(e)}$
and $gexe^{-1}\in X_{\tilde u}= gG_{o(e)}$ by edges of length $1/2$ each. We define 
$f_{\tilde{e},\tilde{v}}:X_{\tilde{e}}\map X_{\tilde{v}}$ and $f_{\tilde{e},\tilde{u}}:X_{\tilde{e}}\map X_{\tilde{u}}$
by setting $f_{\tilde{e},\tilde{v}}(gex)= gex$ and $f_{\tilde{e},\tilde{u}}(gex)= gexe^{-1}$.

We have a natural simplicial map $\pi:X\map \mathcal T$ (more precisely to the first barycentric subdivision of $\mathcal T$).
This map is the coarse analog of the {\em tree of metric spaces} introduced by \cite{BF} (see also \cite{mitra-trees}).
By abuse of terminology we shall refer to this also as a tree of metric spaces or graphs.
We recall some notations and definitions from \cite{mitra-trees} and collect some basic properties.

$(1)$ We note that $X_{u}=\pi^{-1}(u)$ and $X_{\tilde e}=\pi^{-1}( e)$ for all $u\in V(\mathcal T)$ and $e\in E(\mathcal T)$.
For all $u\in V(\mathcal T)$ the intrinsic path metric of $X_u$ will be denoted by $d_u$.
Similarly, we use $d_e$ for the intrinsic path metric on $X_e$. It follows that with these intrinsic metrics the metric spaces
$X_e$, $X_u$ are isometric to the Cayley graphs $\Gamma(G_e,S_e)$ and $\Gamma(G_u,S_u)$ respectively. Therefore, if all
the vertex and edge groups are Gromov hyperbolic then the vertex and edge spaces of $X$ are uniformly hyperbolic metric spaces.

$(2)$ {\bf Quasi-isometric lifts of geodesics:} Suppose $u, v\in \mathcal T$ and let $[u,v]$ denote the geodesic in $\mathcal T$ joining them.
A $K$-QI section of $\pi$ over $[u,v]$ or a $K$-QI lift of $[u,v]$ (in $X$) is a set theoretic section $s:[u,v]\map X$
of $\pi$ which is also a $K$-QI embedding. In general, we are only interested in defining these sections over the vertices in
$[u,v]$. 

$(3)$ {\bf Hallways flare condition}: {\em
We will say that $\pi: X\map \mathcal T$ satisfies the hallways flare condition if for all $K\geq 1$ there are numbers
$\lambda_K>1, M_K\geq 1, n_K\geq 1$ such that given a geodesic $\alpha:[-n_K,n_K]\map \mathcal T$ and two $K$-QI lifts
$\alpha_1,\alpha_2$ of $\alpha$, if $d_{\alpha(0)}(\alpha_1(0),\alpha_2(0))\geq M_K$ then {\small
$$max\{d_{\alpha(n_K)}(\alpha_1(n_K), \alpha_2(n_K)), d_{\alpha(-n_K)}(\alpha_1(-n_K), \alpha_2(-n_K))\}\geq \lambda_K d_{\alpha(0)}(\alpha_1(0),\alpha_2(0)).$$}}
$(4)$ {\em {\bf Graphs of groups with QI embedded conditions.}
Suppose $(\mathcal G, Y)$ is a graph of groups such that each vertex and edge group is finitely generated. We say that
it satisfies the QI embedded condition if all the inclusion maps of the edge groups into the vertex groups
are quasi-isometric embeddings with respect any choice of finite generating sets for the vertex and edge groups.}

It is clear that if $(\mathcal G, Y)$ is a graph of groups with QI embedded condition then all the maps $f_{e,u}:X_e\map X_u$
are uniform QI embeddings.

\begin{lemma}
There is a naturally defined proper and cocompact action of $G$ on $X$ such that the map $\pi:X\map \mathcal T$
is $G$-equivariant.
\end{lemma}
{\em Proof:} We note that $X$ is obtained from the disjoint union of the cosets of the vertex
and edge groups of $(\mathcal G, Y)$. The group $G$ has a natural action on this disjoint union.
It is also easy to check that under this action adjacent vertices of $X$ go to adjacent vertices. Thus we have a simplicial
$G$-action on $X$. Clearly the natural map $\pi:X\map \mathcal T$ is $G$-equivariant.
To show that the action is proper it is enough to show that the vertex stabilizers are uniformly finite. However,
if a point $x\in gG_v$ is fixed by an element $h\in G$ then $h$ fixes $gG_v\in V(\mathcal T)$.
However, stabilizers of $gG_v$ is simply $gG_vg^{-1}$ and the action of $gG_vg^{-1}$ on $gG_v\subset X$ is fixed point free.
Hence, the $G$-action on $X$ is fixed point free.

That the $G$-action is cocompact on $X$ follows from the fact that the $G$-actions on $V(\mathcal T)$ and
$E(\mathcal T)$ are cofinite. $\Box$

Fix a vertex $v_0\in Y$ and the vertex $G_{v_0}\in V(\mathcal T)$. Look at the corresponding vertex space $G_{v_0}\subset X$ and
let $x_0$ denote $1\in G_{v_0}$. Let $\Theta:G\map X$ denote the orbit map $g\mapsto gx_0$. By Milnor-Schwarz lemma this orbit
map is a quasi-isometry since the $G$-action is proper and cocompact by the above lemma.
\begin{lemma}\label{model-map}
There is a constant $D_0$ such that for all vertex space $gG_v\subset X$ we have $Hd(\Theta(gG_v), gG_v)\leq D_0$.

It follows that for any $g.x\in gG_v\subset X$ we have $g.x\in \Theta^{-1}(B(gx,D_0))$.
\end{lemma}
$Proof:$ 
For proving the lemma let $\gamma_v$ be a geodesic in $X$ joining $x_0$ to the identity element of $G_v$. Then for all $x\in G_v$,
$gx\gamma_v$ is a path joining $gxx_0$ and $gx\in gG_v$. Hence one can choose $D_0$ to be the maximum of the lengths of
$\gamma_v$'s, $v\in V(Y)$. $\Box$ 

The following corollary is an immediate consequence of the above two lemmas.

\begin{cor}
The vertex spaces and edge spaces of $X$ are uniformly properly embedded in $X$.
\end{cor}

{\bf Notation:} We shall use $i_w:X_w\map X$ denote the canonical inclusion of the vertex and edge spaces of $X$ into $X$.
Let $\tilde{v}=gG_v\in V(\mathcal T)$. It follows from the above corollary that $\Theta $ induces a coarsely well-defined quasi-isometry
from $gG_v\subset G$ to $X_{\tilde{v}}$. Namely, we can send any $x\in gG_v$ to a point $y$ of $X_{\tilde{v}}$
such that $d_X(\Theta(x),y)\leq D_0$. where $D_0$ is as in the above corollary. 
We shall denote this by $\Theta_{g,v}:gG_v\map X_{\tilde{v}}$.

\section{The main theorem}

For the rest of the paper we shall assume that $G$ is a hyperbolic group which admits a graph of groups decomposition
$(\mathcal G, Y)$ with the QI embedded condition where all the vertex and edge groups are hyperbolic.
Let $\mathcal T$ be the Bass-Serre tree of this graph of groups.

We aim to show that in $G$ the family of subgroups $\{G_v: v\in V(\mathcal T)\}$ satisfies
the limit set intersection property:

\begin{theorem}\label{main-thm}
Suppose a hyperbolic group $G$ admits a decomposition into a graph of hyperbolic groups $(\mathcal G,Y)$ with
quasi-isometrically embedded condition and suppose $\mathcal T$ is the correspoding Bass-Serre tree.
Then for all $w_1, w_2\in V(\mathcal T)$ we have $\Lambda(G_{w_1})\cap \Lambda(G_{w_2})= \Lambda(G_{w_1}\cap G_{w_2})$.
\end{theorem}

The idea of the proof is to pass to the tree of space $\pi:X\map \mathcal T$ using the orbit map $\Theta:G\map X$
defined in the previous section and then use the techniques of \cite{mitra-trees}. 
The following theorem is an important ingredient of the proof.

\begin{theorem} (\cite{mitra-trees})\label{mit-thm}
The inclusion maps $i_w: X_w\map X$ admit CT maps $\partial i_w: \partial X_w \map \partial X$
for all $w\in V(\mathcal T)$.
\end{theorem}

Recall that if $u,v\in V(\mathcal T)$ are connected by an edge $e$ there are natural maps $f_{e,u}:X_e\map X_u$ and
$f_{e,v}:X_e\map X_v$. We know that these maps are uniform QI embeddings. We assume that they are  all $K$-QI embeddings for some $K>1$.
They induce embeddings $\partial f_{e,u}:\partial X_e\map \partial X_u$ and $\partial f_{e,v}:\partial X_e\map \partial X_v$
by Lemma \ref{bdry-top}. Therefore, we get
partially defined maps from $\partial X_u$ to $\partial X_v$ with domain $Im(\partial f_{e,u})$. Let us denote this
by $\psi_{u,v}:\partial X_u\map \partial X_v$. By definition for all $x\in X_e$ we have
$\psi_{u,v}(f_{e,u}(x))=f_{e,v}(x)$.

\begin{defn}
$(1)$ If $\xi\in \partial X_u$ is in the domain of $\psi_{u,v}$ and $\psi_{u,v}(\xi)=\eta$ then we say that
$\eta$ is a flow of $\xi$ and that $\xi$ can be flowed to $\partial X_v$.

$(2)$ Suppose $w_0\neq w_n\in V(\mathcal T)$ and $w_0,w_1,\cdots, w_n$ are consecutive vertices of the geodesic
$[w_0,w_n]\subset \mathcal T$. We say that a point $\xi\in \partial X_{w_0}$ can be flowed to $\partial X_{w_n}$ if
there are $\xi_i\in \partial X_{w_i}$, $0\leq i\leq n$ where $\xi_0=\xi$ such that $\xi_{i+1}=\psi_{w_i,w_{i+1}}(\xi_i)$,
$0\leq i \leq n-1$.
In this case, $\xi_n$ is called the flow of $\xi_0$ in $X_{w_n}$. 
\end{defn}

Since the maps $\psi_{u,v}$ are injective on the their domains for all $u\neq v\in V(\mathcal T)$ and $\xi\in \partial X_u$ the flow of
$\xi$ in $\partial X_v$ is unique if it exists.

\begin{lemma} \label{main-lemma}
Suppose $w_1,w_2\in V(\mathcal T)$ and $\xi_i\in \partial X_{w_i}$, $i=1,2$ such that $\xi_2$ is a flow of $\xi_1$.
Let $\alpha_i$ be a geodesic in the vertex space $X_{w_i}$ such that $\alpha_i(\infty)=\xi_i$, $i=1,2$.
Then $Hd(\alpha_1, \alpha_2)< \infty$.
\end{lemma}

$Proof:$ It is enough to check it when $w_1,w_2$ are adjacent vertices. Suppose $e$ is the edge connecting $w_1,w_2$.
The lemma follows from the stability of quasi-geodesics in the hyperbolic space $X_{w_i}$'s and the fact that
every point of $X_e$ is at distance $1/2$ from $X_{w_i}$, $i=1,2$. $\Box$

\begin{cor}\label{main-cor}
Under the CT maps $\partial i_{w_j}: \partial X_{w_j}\map  \partial X$, $j=1,2$ the points $\xi_1,\xi_2$
go the same point of $\partial X$, i.e. $\partial i_{w_1}(\xi_1)=\partial i_{w_2}(\xi_2)$.
\end{cor}

\begin{lemma}\label{bdd-flow}
Let $w_1,w_2\in \mathcal T$ and suppose they are joined by an edge $e$. Suppose $\xi_1\in \partial X_{w_1}$ can not be
flowed to $\partial X_{w_2}$. Let $\alpha\subset X_{w_1}$ be a geodesic ray such that $\alpha(\infty)=\xi_1$.
Then for all $D>0$ the set $N_D(\alpha)\cap f_{e,w_1}(X_e)$ is bounded.
\end{lemma}

$Proof:$ If $N_D(\alpha)\cap f_{e,w_1}(X_e)$ is not bounded for some $D>0$ then $\xi_1$ is in the limit set of 
$f_{e,w_1}(X_e)$ and so $\xi_1$ can be flowed to $\partial X_{w_2}$ by Lemma \ref{CT-limset}.
This contradiction proves the lemma. $\Box$

Now we briefly recall the {\em ladder construction} of Mitra which was crucial for the proof of 
the main theorem of \cite{mitra-trees}. We shall need it for the proof of Theorem \ref{main-thm}.

{\bf Mitra's Ladder} $B(\lambda)$.

Fix $D_0,D_1>0$. Let $v\in V(\mathcal T)$ and $\lambda$ be a finite geodesic segment of $X_v$. We shall define the set
$B(\lambda)$ to be a union of vertex space geodesics $\lambda_w\subset X_w$ where $w$ is in a subtree $T_1$ of $\mathcal T$ containing $v$.
The construction is inductive. Inductively one constructs the $n$-sphere $S_{T_1}(v,n)$ of $T_1$ centred at $v$ and the corresponding
$\lambda_w$'s, $w\in S(v,n)$.

$S_{T_1}(v,1)$: There are only
finitely many edges $e$ incident on $v$ such that $N_{D_0}(\lambda)\cap f_{e,v}(X_e)\neq \emptyset$.
Then $S(v,1)$ is the set of terminal points of all the edges $e$ that start at $v$ such that the diameter of 
$N_{D_0}(\lambda)\cap f_{e,v}(X_e)$ is at least $D_1$. In this case, for each edge $e$ connecting $v$ to say $v_1\in S(v,1)$,
we choose two points, say $x,y\in N_{D_0}(\lambda)\cap f_{e,v}(X_e)$ such that $d_v(x,y)$ is maximum. Then we
choose $x_1, y_1\in X_{v_1}$ such that $d(x,x_1)=1$ and $d(y,y_1)=1$ and define $\lambda_{v_1}$ to be a geodesic
in $X_{v_1}$ joining $x_1,y_1$.

$S_{T_1}(v,n+1)$ from $S_{T_1}(v,n)$: Suppose $w_1\in S(v,n)$. Then a vertex $w_2$ adjacent to $w_1$ with $d_T(v,w_2)=n+1$
belongs to $S(v,n+1)$ if the diameter of $N_{D_0}(\lambda_{w_1})\cap f_{e,w_1}(X_e)$ is at least $D_1$,
where $e$ is the edge connecting $w_1,w_2$, has diameter at least $D_1$ in $X_{w_1}$.
To define $\lambda_{w_2}$ one chooses two points $x,y\in N_{D_0}(\lambda_{w_1})\cap f_{e,w_1}(X_e)$ such that $d_{w_1}(x,y)$
is maximum, then let $x_1,y_1\in X_{w_2}$ be such that $d(x,x_1)=1$ and $d(y,y_1)=1$ and define $\lambda_{w_2}$ to be a geodesic
in $X_{w_2}$ joining $x_1,y_1$.

\begin{theorem}(Mitra \cite{mitra-trees})\label{mit-ladd}
There are constants $D_0>0, D_1>0$ and $C>0$ depending on the defining parameters of the tree of metric spaces
$\pi:X\map \mathcal T$ such that the following holds:

For any $v\in V(\mathcal T)$ and a geodesic segment $\lambda \subset X_v$ the corresponding ladder $B(\lambda)$ is a
$C$-quasi-convex subset of $X$. 
\end{theorem}

To prove this theorem, Mitra defines a coarse Lipschitz retraction map $P: X\map B(\lambda)$ which we now recall.
For the proof of how this works one is referred to \cite{mitra-trees}. {\em However, we shall subsequently assume that
appropriate choices of $D_0,D_1$ are made in our context so that all the ladders are uniformly quasiconvex subsets
of $X$.}

{\bf Coarsely Lipschitz retraction on the ladders}

Suppose $\lambda \subset X_v$ is a geodesic. Let $T_1=\pi(B(\lambda))$.
For each $w\in T_1$, $\lambda_w=X_w\cap B(\lambda)$ is a geodesic in $X_w$. We know that there is a coarsely well
defined nearest point projection $P_w: X_w\map \lambda_w$. (See Proposition $3.11$ in Chapter $III.\Gamma$ of \cite{bridson-haefliger}.)
Now for each $x\in X_w$, $w\in T_1$ define
$P(x)=P_w(x)$. If $x\in X_w$ and $w\not \in T_1$ then connect $w$ to $T_1$ by a geodesic in $\mathcal T$.
Since $\mathcal T$ is a tree there is a unique such geodesic. Let $w_1\in T_1$ be the end point of this geodesic
and let $e$ be the edge on this geodesic incident on $w_1$ going out of $T_1$. Mitra proved that in this case the projection of
$f_{e,w_1}(X_e)$ on $\lambda_{w_1}$ is uniformly small. It follows by careful choice of $D_0,D_1$. (See Lemma $3.1$ in \cite{mitra-trees}.)
Choose a point $x_{w_1}$ on this projection. Define $P(x)=x_{w_1}$.

\begin{theorem}(\cite{mitra-trees})\label{mit-proj}
The map $P:X\map B(\lambda)$ is a coarsely Lipschitz retraction. 
\end{theorem}
In other words, it is a retraction and there are constants $A,B$ such that $d(P(x),P(y))\leq Ad(x,y)+B$ for all $x,y\in X$.

Using the above theorems of Mitra now we shall prove the converse of Corollary \ref{main-cor}.
This is the last ingredient for the proof of Theorem \ref{main-thm}.

\begin{prop}\label{main-prop}
Suppose $v\neq w\in \mathcal T$ and there are points $\xi_v\in \partial X_v$ and
$\xi_w \in \partial X_w$ which map to the same point $\xi\in \partial X$ under the CT maps
$\partial {X}_v\rightarrow \partial {X}$ and $\partial {X}_w\rightarrow \partial {X}$ respectively. Then $\xi_v$ can
be flowed to $\partial X_w$.
\end{prop}

$Proof:$ Using Lemma \ref{main-lemma} we can assume that the point $v$ is such that $\xi_v$ can not
further be flowed along $vw$ and similarly $\xi_w$ can not be flowed in the direction of $wv$ where $v\neq w$.
Let $\alpha: [0,\infty) \rightarrow X_v$ and $\beta: [0,\infty)\rightarrow X_w$ be geodesic rays in $X_v,X_w$ respectively
such that $\alpha(\infty)=\xi_v$ and $\beta(\infty)=\xi_w$.
Let $e_v, e_w$ be the first edges from the points $v, w$  along the direction of $vw$ and $wv$
respectively. Then $N_{D_0}(\alpha)\cap f_{e_v,v}(X_{e_v})\subset X_v$ and $N_{D_0}(\beta)\cap f_{e_w,w}(X_{e_w})\subset X_w$
are both bunded sets by Lemma \ref{bdd-flow} where $D_0$ is as in Theorem \ref{mit-ladd}.

For all $n\in \NN$ let $\alpha_n:= \alpha|_{[0, n]}$ and
$\beta_n:=\beta|_{[0, n]}$ respectively. The ladders $B(\alpha_n), B(\beta_n)$ are uniformly quasi-convex
subsets of $X$ by Theorem \ref{mit-proj}. Hence there
are uniform ambient quasi-geodesics of $X$ in these ladders joining $\alpha(0),\alpha(n)$ and $\beta(0), \beta(n)$
respectively. Choose one such for each one of them and let us call them $\gamma_n$ and $\gamma^{'}_n$
respectively. Now, since $\alpha$ and $\beta$ limit on the same point $\xi\in \partial X$, 
by Lemma \ref{bdry-lemma}(3) there is a uniform constant $D$ such that for a subsequence $\{n_k\}$ of natural numbers
there are points $x_{n_k}\in \gamma_{n_k}$ and $y_{n_k}\in \gamma^{'}_{n_k}$ such that $d(x_{n_k},y_{n_k})\leq D$
and $\lim_{k\map\infty} x_{n_k}=\lim_{k\map \infty}y_{n_k}=\xi\in\partial X$.

Let $v_1, w_1$ be the vertices on the geodesic $[v,w]\subset \mathcal T$ adjacent to $v,w$ respectively.
Let $A_{v_k}:=B(\alpha_{n_k})\cap X_{v_1}$ and $A_{w_k}:=B(\beta_{n_k})\cap X_{w_1}$ respectively. 

If we remove the edge space $X_{e_v}$ from $X$ then 
the remaining space has two components- one containing $X_v$ and the other containing $X_w$.
Call them $Y_1,Y_2$ respectively. We note that since the diameter of $A_{v_k}$ is uniformly bounded, if at all nonempty,
the portion of $\gamma_{n_k}$ contained in $Y_2$, if at all it travels into $Y_2$, is uniformly bounded. This implies that
the portion of $\gamma_{n_k}$ joining $\alpha(0)$ and $A_{v_k}$ is uniformly small if $A_{v_k}\neq \emptyset$.

Hence, there are infinitely many $k\in \NN$ such that $x_{n_k}\in Y_1$ and $y_{n_k}\in Y_2$. Since we are dealing with
a tree of spaces and $d(x_{n_k}, y_{n_k})\leq D$ for all $k\in \NN$, this implies there are points
$z_k\in f_{e_v}(X_{e_v})$ such that $d(x_{n_k},z_k)\leq D$ for all $k\in \NN$.
Thus $\xi_v$ can be flowed to $\partial X_{v_1}$ by
Lemma \ref{CT-limset}. This contradiction proves the proposition. $\Box$

{\em Proof of Theorem \ref{main-thm}}:

Suppose $w_i=g_iG_{v_i}$, $i=1,2$, for some $g_1,g_2\in G$ and $v_1,v_2\in V(Y)$. This implies
$G_{w_i}=g_iG_{v_i}g^{-1}_i$, $i=1,2$. Also, $\Lambda(g_iG_{v_i}g^{-1}_i)=\Lambda(g_iG_{v_i})$ by Lemma \ref{int-lemma}(1).
Hence we need to show that $\Lambda(g_1G_{v_1})\cap \Lambda(g_2G_{v_2})= \Lambda(g_1G_{v_1}g^{-1}_1\cap g_2G_{v_2}g^{-1}_2)$.
Using Lemma \ref{int-lemma}(2), therefore, it is enough to show that
$$\Lambda(G_{v_1})\cap \Lambda(gG_{v_2})= \Lambda(G_{v_1}\cap gG_{v_2}g^{-1})\,\, \mbox{for all}\,\, v_1,v_2\in V(Y), g\in G.$$

Clearly, $ \Lambda(G_{v_1}\cap gG_{v_2}g^{-1})\subset \Lambda(G_{v_1})\cap \Lambda(gG_{v_2})$. Thus we need to show
that $\Lambda(G_{v_1})\cap \Lambda(gG_{v_2})\subset \Lambda(G_{v_1}\cap gG_{v_2}g^{-1})$.

Given an element $\xi\in \Lambda(G_{v_1})\cap \Lambda(gG_{v_2})$
there are $\xi_1\in \partial G_{v_1}$ and $\xi_2\in \partial (gG_{v_2})$ both of which map to $\xi$ under the CT maps
$\partial G_{v_1}\map \partial G$ and $\partial (gG_{v_2})\map \partial G$ by Lemma \ref{CT-limset}.
Now, we have a quasi-isometry $\Theta:\Gamma(G,S)\map X$. By Lemma \ref{model-map} each coset of any
vertex group in $G$ is mapped uniformly Hausdorff close to the same coset in $X$. Hence $\Theta$ induces uniform quasi-isometries
$\Theta_{g,v}$ from $gG_v\subset \Gamma(G,S)$ to $gG_v\subset X$ for all $g\in G$, $v\in Y$. For avoiding confusion
let us denote the subset $G_{v_1}\subset X$ by $X_{w_1}$ and $gG_{v_2}\subset X$ by $X_{w_2}$. It follows that
$\partial \Theta_{1,v_1}(\xi_1)\in \partial X_{w_1}$ and $\partial\Theta_{g,v_2}(\xi_2)\in \partial X_{w_2}$ are mapped
to the same element of $\partial X$ under the CT maps $\partial X_{w_i}\map \partial X$, $i=1,2$.
 Hence by Proposition \ref{main-prop} $\partial \Theta_{1,v_1}(\xi_1)$ can be flowed to, say $\xi^{'}_2\in
\partial X_{w_2}$. By Lemma \ref{main-lemma} the image of $\xi^{'}_2$ and $\Theta_{g,v_2}(\xi_2)$ under the
CT map $\partial X_{w_2}\map \partial X$ are the same. Hence, we can replace $\xi_2$ by $(\partial \Theta_{g,v_2})^{-1}(\xi_2)$
and assume that $\Theta_{1,v_1}(\xi_1)$ flows to $\Theta_{g,v_2}(\xi_2)$. 

Then by Lemma \ref{main-lemma} for any geodesic rays $\alpha_i \subset X_{w_i}$ with
 $\alpha_i(\infty)= \partial \Theta_{x_i,v_i}(\xi_i)$ for $i=1,2$ where $x_1=1$ and $x_2=g$ we have  $Hd( \alpha_1, \alpha_2)< \infty$.
Pulling back these geodesics by $\Theta_{1,v_1}$ and $\Theta_{g,v_2}$ we get uniform quasi-geodesic rays,
say $\beta_1\subset G_{v_1}$ and $\beta_2\subset  gG_{v_2}$ such that $\beta_i(\infty)=\xi_i$, $i=1,2$ and $Hd(\beta_1,\beta_2)<\infty$.

Now let $p_i=\beta_1(i)$ and $q_i\in \beta_2$, $i\in \mathbb N$  be such that $d(p_i,q_i)\leq D$ where $Hd(\beta_1,\beta_2)=D$.
Join $p_i$ to $q_i$ by a geodesic in $\Gamma(G,S)$. Suppose $w_i$ is the word labeling this geodesic. Since there are only finitely
many possibilities for such words, there is a constant subsequence $\{w_{n_k}\}$ of $\{w_n\}$. Let $h_k=p^{-1}_{n_1}p_{n_k}$
and $h^{'}_k=q^{-1}_{n_1}q_{n_k}$. Let $x$ be the group element represented by $w_{n_k}$.
Then we have $p_{n_1}.h_k.x=p_{n_1}.x.h^{'}_k$ or $h^{'}_k=xh_kx^{-1}$. Since $h^{'}_k$ connects two elements of $gG_{v_2}$,
it is in $G_{v_2}$. Hence $h_k\in G_{v_1}\cap xG_{v_2}x^{-1}$. Thus 
$p_{n_1}h_kp^{-1}_{n_1}\in p_{n_1}G_{v_1}p^{-1}_{n_1} \cap p_{n_1}xG_{v_2}(p_{n_1}x)^{-1}=G_{v_1}\cap gG_{v_2}g^{-1}$.
Finally, since $d(p_{n_1}h_kp^{-1}_{n_1}, p_{n_1}h_k)=d(p_{n_1}h_kp^{-1}_{n_1}, p_{n_k})=d(1,p_{n_1})$ for all $k\in \NN$,
$\lim_{n\map \infty} p_{n_1}h_kp^{-1}_{n_1}= \lim_{n \map \infty} p_{n_k}=\xi_1$. This completes the proof. $\Box$

The following corollary has been pointed out by Mahan Mj. We use the same notations as in the main theorem.

\begin{cor}
If $H_i\subset G_{w_i}$, $i=1,2$ are two quasiconvex subgroups then $\Lambda(H_1)\cap \Lambda(H_2)=\Lambda(H_1\cap H_2)$.
\end{cor}

$Proof:$ Assume that $w_i=g_iG_{v_i}$, $i=1,2$. Then $G_{w_i}=g_iG_{v_i}g^{-1}_i$. Let $K_i=g^{-1}_iH_ig_i< G_{v_i}$.
We may construct a new finite graph starting from $Y$ by adding two vertices
$u_1, u_2$ where $u_i$ is connected to $v_i$ by an edge $e_i$, $i=1,2$. Let us call this graph $Y_1$. Define a new
graph of groups $(\mathcal G_1,Y_1)$ by keeping the definition same on $Y$ and setting $G_{u_i}=G_{e_i}=K_i$, $i=1,2$
and defining $\phi_{e_i,u_i}=1_{K_i}$ and $\phi_{e_i,v_i}$ to be the inclusion map $K_i\subset G_{v_i}$.
This produces a new graph of groups
with QI embedded condition and with fundamental group isomorphic to $G$. Suppose the Bass-Serre tree of the new graph
of groups is $\mathcal T_1$.

Now we can apply Theorem \ref{main-thm} to $G_{w^{'}_i}$, $i=1,2$ where $w^{'}_i=g_iK_i\in \mathcal T_1$, $i=1,2$
to finish the proof. $\Box$

\begin{example}
We now give an example where intersection of limit sets is not equal to the limit set of the intersection.
Suppose $G$ is a hyperbolic group with an infinite normal subgroup $H$ such that $G/H$ is not torsion. Let $g\in G$
be such that its image in $G/H$ is an element of infinite order. Let $K=<g>$. Then $H\cap K=(1)$ whence
$\Lambda(H\cap K)=\emptyset$. However, $H$ being an infinite normal subgroup of $G$ we have $\Lambda(H)=\partial G$.
Thus $\Lambda(H)\cap \Lambda(K)=\Lambda(K)\neq \emptyset$.

\end{example}

We end with a question.
\begin{qn}
If a hyperbolic group $G$ admits a decomposition into a graph of hyperbolic groups with qi embedded condition
and $G_v$ is a vertex group how to describe $\Lambda(G_v)\subset \partial G$?
\end{qn}
It has been pointed out to me by Prof. Ilya Kapovich that the first interesting case where this question should be
considered is a hyperbolic strictly ascending HNN extension of a finitely generated nonabelian free group $F$
$$G=<F, t| t^{-1} wt=\phi(w),\,\forall w\in F>$$ where $\phi:F\map F$ is an injective but not surjective endomorphism
of $F$. One would also like to describe $\partial G$ in this case.

\bibliography{lim-int.bib}
\bibliographystyle{amsalpha}

\end{document}